\documentclass[12pt,a4paper,reqno]{amsart}

\usepackage{anysize}
\usepackage[utf8]{inputenc}
\usepackage{amsmath}
\usepackage{amssymb}
\usepackage{array}

\marginsize{2.5cm}{2.5cm}{3cm}{3cm}

\newtheorem{theorem}{Theorem}[section]

\newtheorem{example}[theorem]{Example}

\begin{document}

\title[Heuristic derivation of Zudilin's supercongruences]{Heuristic derivation of Zudilin's supercongruences \\ for rational Ramanujan series}

\author{Jesús Guillera}
\address{Department of Mathematics, University of Zaragoza, 50009 Zaragoza, SPAIN}
\email{}
\date{}

\begin{abstract}
We derive, using a heuristic method, a $p$-adic mate of bilateral Ramanujan series. It has (among other consequences) the Zudilin's supercongruences for rational Ramanujan series. 
\end{abstract}

\maketitle

\section{Rational Ramanujan series for $\pi^{-m}$}  
At the beginning of the twenty first century we discovered new families of Ramanujan-like series, but of greater degree \cite{Rama-new}, and proved several of them by the WZ (Wilf-Zeilberger) method \cite{Rama-WZ}. \par
We can write the rational Ramanujan-like series as
\[
\sum_{n=0}^{\infty} R(n)=\sum_{n=0}^{\infty} \left( \prod_{i=0}^{2m} \frac{(s_i)_n}{(1)_n} \right) \sum_{k=0}^m a_k n^k z_0^n \, = \frac{\sqrt{(-1)^m \chi}}{\pi^m},
\]
where $z_0$ is a rational such that $z_0 \neq 0$ and $z_0 \neq 1$, the parameters $a_0, a_1,...,a_m$ are positive rationals, and $\chi$ the discriminant of a certain quadratic field (imaginary or real), which is an integer. In case that $|z_0|>1$ we understand the series as its analytic continuation. Below, we show an example
\[
\frac{1}{2048} \sum_{n=0}^{\infty}  \frac{\left(\frac{1}{2}\right)_n^7 \left(\frac{1}{4}\right)_n \left(\frac{3}{4}\right)_n}{(1)_n^9} (43680n^4+20632n^3+4340n^2+466n+21) \left( \frac{1}{2^{12}} \right)^n =\frac{1}{\pi^4},
\]
conjectured by Jim Cullen, and recently proved by Kam Cheong Au, using the WZ method \cite{Kam-Cheong-Au}.

\section{Bilateral Ramanujan series}  
We define the function $f: \mathbb{C} \longrightarrow \mathbb{C}$ in the following way:
\[
f(x)=e^{-i \pi x} \prod_{s_k} \frac{\cos \pi x - \cos \pi s_k}{1-\cos \pi s_k} \sum_{n \in \mathbb{Z}} R(n+x).
\]
Then, there exists coefficients $\alpha_k$ and $\beta_k$ (which we conjecture are rational) such that $f(x)=F(x)$, where
\[
F(x)=\frac{\sqrt{(-1)^m \chi}}{\pi^m} \, \left( 1 - \sum_{k=1}^{m} \left(\alpha_k (\cos 2 \pi k x -1) + \beta_k \sin 2 \pi k x \right) \right),
\]
is the Fourier expansion of $f(x)$.
\begin{proof}
The function $f(x)$ is periodic of period $x=1$ because the product over $s_k$ is periodic as each $s_k=s$, has a companion $s_k=1-s$, and the sum over $\mathbb{Z}$ is clearly periodic as well. In addition $f(x)$ is holomorphic because the zero of $\cos \pi x - \cos \pi s_k$ at $x=-s_k$ cancels the pole of $(s_k)_{n+x}$ at $x=-s_k$, and as $f(x)$ is periodic all the other poles are canceled as well. As $f(x)$ is holomorphic and periodic, it has a Fourier expansion. Finally, we can prove that $f(x)=\mathcal{O}(e^{(2m+1) \pi |{\rm Im}(x)|})$, and therefore the Fourier expansion terminates at $k=m$.
\end{proof} 

\begin{example} \rm
\begin{multline}\nonumber 
\frac18 \sum_{n \in \mathbb{Z}} \frac{\left(\frac12\right)_{n+x}^5}{(1)_{n+x}^5} (20(n+x)^2+8(n+x)+1) \left(\frac{-1}{4}\right)^{n+x} \\ =
e^{i \pi x} \frac{1-\frac12 (\cos 2 \pi x -1) + \frac12 (\cos 4 \pi x -1)}{\pi^2 \cos^5\pi x}.
\end{multline}
\end{example}

\begin{example} \rm
\begin{multline*}
\frac{1}{384} \sum_{n \in \mathbb{Z}} \frac{\left(\frac12\right)_{n+x}\left(\frac13\right)_{n+x}\left(\frac23\right)_{n+x}\left(\frac16\right)_{n+x}\left(\frac56\right)_{n+x}}{(1)_{n+x}^5}  \left( - \frac{3^6}{4^6} \right)^{n+x} \\ \times \left( 1930(n+x)^2+549(n+x)+45 \right) \\ = 
e^{i \pi x} \frac{3-14 (\cos 2\pi x - 1) + 6 (\cos 4\pi x - 1)}{\pi^2 \cos \pi x \, (4\cos^2 \pi x-1)(4\cos^2 \pi x -3)}.
\end{multline*}
\end{example}

\begin{example} \rm
\begin{multline*}
\frac{1}{32} \sum_{n \in \mathbb{Z}} \frac{\left(\frac12\right)_{n+x}^3\left(\frac14\right)_{n+x}\left(\frac34\right)_{n+x}}{(1)_{n+x}^5}  \left( \frac{1}{16} \right)^{n+x} (120(n+x)^2+34(n+x)+3) = \\
e^{i \pi x} \, \frac{1-\frac72 (\cos 2 \pi x -1) + \frac32 (\cos 4 \pi x - 1) + \left( \frac12 \sin 2\pi x - \frac12 \sin 4 \pi x \right) i}{\pi^2 \cos^3 \pi x \, (2\cos^2 \pi x -1)}.
\end{multline*}
\end{example}

\begin{example} \rm 
\begin{multline*}
\frac16 \sum_{n \in \mathbb{Z}} \frac{(\frac12)_{n+x}^3(\frac13)_{n+x}(\frac23)_{n+x}}{(1)_{n+x}^5} \left[ 28(n+x)^2+18(n+x)+3 \right] \,  (-27)^{n+x} \\ 
= e^{i \pi x} \frac{3 + (\cos 2\pi x - 1) + \frac34 (\cos 4\pi x - 1)}{\pi^2 \cos^3 \pi x (4\cos^2 \pi x - 1)}.
\end{multline*}
\end{example}

\section{Series to the right and to the left}
The series to the right hand side is
\[
\sum_{n=0}^{\infty} R(n+x) = \sum_{n=0}^{\infty} \left( \prod_{i=0}^{2m} \frac{(s_i)_{n+x}}{(1)_{n+x}} \right) \sum_{k=0}^m a_k (n+x)^k \, z_0^{n+x},
\]
extended by analytic continuation to all $z_0$ different from $0$ and $1$, and
the series to the left hand side is
\begin{align*}
\sum_{n=1}^{\infty} R(-n+x) &= \sum_{n=1}^{\infty} \left( \prod_{i=0}^{2m} \frac{(s_i)_{-n+x}}{(1)_{-n+x}} \right) \sum_{k=0}^m a_k (-n+x)^k \, z_0^{-n+x}
\\ &= 
x^{2m+1} \sum_{n=1}^{\infty} \left( \prod_{i=0}^{2m} \frac{(1)_{n-x}}{(s_i)_{n-x}} \right) \sum_{k=0}^m a_k (-n+x)^{k-2m-1} \, z_0^{-n+x}
\\ &=
x^{2m+1} z_0^{x} \left( \prod_{i=0}^{2m} \frac{(s_i)_{x}}{(1)_{x}} \right) \sum_{n=1}^{\infty} \left( \prod_{i=0}^{2m} \frac{(1-x)_{n}}{(s_i-x)_{n}} \right) \sum_{k=0}^m a_k (-n+x)^{k-2m-1} \, z_0^{-n},
\end{align*}
extended by analytic continuation to all $z_0$ different from $0$ and $1$. We see that
\begin{multline}\nonumber
\sum_{n=0}^{\infty} R(n) - \sum_{n=0}^{\infty} R(n+x) =  \frac{\sqrt{(-1)^m \chi}}{\pi^m}
\\ 
- e^{i \pi x} \prod_{s_k} \frac{1-\cos \pi s_k}{\cos \pi x - \cos \pi s_k} \frac{\sqrt{(-1)^m \chi}}{\pi^m} \left(1 - \sum_{k=1}^{m} \left(\alpha_k (\cos 2 \pi k x -1)  + \beta_k \sin 2 \pi k x \right) \right) \\ + (A + B x + C x^2 + \cdots)x^{2m+1}, \quad |x|<1,
\end{multline}
where $(A + B x + C x^2 + \cdots)x^{2m+1}$ is the development of the series to the left hand side at $x=0$, that is 
\[
z_0^{x} \left( \prod_{i=0}^{2m} \frac{(s_i)_{x}}{(1)_{x}}  \right) \sum_{n=1}^{\infty} \left( \prod_{i=0}^{2m} \frac{(1-x)_{n}}{(s_i-x)_{n}} \right) \sum_{k=0}^m a_k (-n+x)^{k-2m-1} \, z_0^{-n} = A + B x + C x^2 + \cdots.
\]

\section{Heuristic derivation of a $p$-adic mate}
Let 
\[
S(N) = \sum_{n=0}^ {\infty} R(n) - \sum_{n=0}^{\infty} R(n+N) = \sum_{n=0}^{N-1} R(n).
\]
As in a Ramanujan-like series each $s_k<1/2$ has a companion $1-s_k$, we notice that 
\[
e^{i \pi x} \prod_{s_k} \frac{1-\cos \pi s_k}{\cos \pi x - \cos \pi s_k}
= e^{i \pi x} \prod_{s_k=\frac12} \frac{1}{\cos \pi x} \prod_{s_k<\frac12} \frac{1-\cos^2 \pi s_k}{\cos^2 \pi x - \cos^2 \pi s_k}
\]
tends to $1$ as $x \to N$ because there is an odd number of factors when $s_k=1/2$. Hence for $x \to N$, we formally have 
\begin{multline}\nonumber
S(x)=\sum_{n=0}^{\infty} R(n) - \sum_{n=0}^{\infty} R(n+x) =  
\\ 
\frac{\sqrt{(-1)^m \chi}}{\pi^m} \left(\sum_{k=1}^{m} \left(\alpha_k (\cos 2 \pi k x -1)  + \beta_k \sin 2 \pi k x \right)\right) 
\\ + (A + B x + C x^2 + \cdots)x^{2m+1}.
\end{multline}
Let
\[
G(x)=\frac{\sqrt{(-1)^m \chi}}{\pi^m} \left(\sum_{k=1}^{m} \left(\alpha_k (\cos 2 \pi k x -1)  + \beta_k \sin 2 \pi k x \right)\right).
\]
For obtaining the $p$-adic analogues $G_p(xp)$ and $G_p(x)$ of $G(xp)$ and $G(x)$, we develop $G(xp)$ and $G(x)$ in powers of $x$. Then, replace the powers of $\pi$ using values of Dirichlet $L$-functions, and the $L$-functions with the corresponding $p$-adic $L$-functions. Finally, we observe that due to the properties of the $L_p$-functions we have: For positive $\chi$ the even powers of $\pi$ turn to $0$, and for negative $\chi$, are the odd powers of $\pi$ that turn to $0$.  This occurs also for negative powers of $\pi$ (regularization). After making the replacements, we see that 
\[
\lim_{x \to \nu} \frac{G_p(xp)}{G_p(x)} = p^m.
\]
For $x=\nu$, where $\nu=1,2,3,\dots$, we see that
\[
z_0^{\nu} \left( \prod_{i=0}^{2m} \frac{(s_i)_{\nu}}{(1)_{\nu}} \right) \sum_{n=1}^{\infty} \left( \prod_{i=0}^{2m} \frac{(1-\nu)_{n}}{(s_i-\nu)_{n}} \right) \sum_{k=0}^m a_k (-n+\nu)^{k-2m-1} \, z_0^{-n} =  A'+B'\nu+C'\nu^2+\cdots,
\]
where 
\[
A'= z_0^{\nu} \left( \prod_{i=0}^{2m} \frac{(s_i)_{\nu}}{(1)_{\nu}} \right) A, \quad 
B'= z_0^{\nu} \left( \prod_{i=0}^{2m} \frac{(s_i)_{\nu}}{(1)_{\nu}} \right) B, \dots.
\]
On the other hand, we see that
\[
S(\nu)=(A'+B'\nu+C'\nu^2+\cdots) \nu^{2m+1} = z_0^{\nu} \left( \prod_{i=0}^{2m} \frac{(s_i)_{\nu}}{(1)_{\nu}} \right) (A+B\nu+C\nu^2+\cdots) \nu^{2m+1}.
\]

\par To get the $p-$adic mate of $S(x)$ we must divide $S_p(\nu p)$ enter $S_p(\nu)$, taking into account that the contribution of $G(x)$ is $(\chi/p) p^m$, and the contribution of the left hand sum is given by
\begin{multline}\nonumber
\frac{(A_p+B_p\nu p + C_p \nu^2 p^2+\cdots) p^{2m+1} \nu^{2m+1}}{(A+B\nu + C \nu^2+\cdots)_p \nu^{2m+1}} \\
= z_0^{\nu} \left( \prod_{i=0}^{2m} \frac{(s_i)_{\nu}}{(1)_{\nu}} \right)\frac{(A_p+B_p\nu p + C_p \nu^2 p^2+\cdots) p^{2m+1} \nu^{2m+1}}{S_p(\nu)}.
\end{multline}
Observe that $(\chi/p)$ is associated to $S(\nu p)$ and $(\chi/1)=1$ to $S_p(\nu)=S(\nu)$. Finally, we must sum both contributions arriving at the following $p$-adic identity:
\[
S(\nu p) = S(\nu) \left(\frac{\chi}{p}\right) p^m + T(\nu)  (A_p+B_p\nu p + C_p \nu^2 p^2+\cdots) p^{2m+1} \nu^{2m+1},
\]
where $A_p, B_p,C_p \dots$, are the $p$-adic analogues of $A,B,C\dots$, and
\[
T(\nu)= z_0^{\nu} \left( \prod_{i=0}^{2m} \frac{(s_i)_{\nu}}{(1)_{\nu}} \right).
\]
Observe that taking positive integers values of $\nu$ we can eliminate some of the constants $A_q$, $B_q$,.., and obtain a new kind of supercongruences $\pmod{p^{2m+k}}$. For example, eliminating $A_q$ and $B_q$, we obtain supercongruences $\pmod{p^{2m+3}}$ relating $S(p)$, $S(2p)$ and $S(3p)$.
\par We can apply a similar technique of bilateral series and $p$-adic mates to other kind of hypergeometric series, for example to those in \cite{Rama-Orr}.

\section{Extended Zudilin's supercongruences}
The above $p$-adic mate has (among other consequences) a generalization for positive integers $\nu$ of Zudilin's $\nu=1$ supercongruences \cite{Zudilin-supercong} and \cite{Rama-divergent-supercong}, namely
\[
S(\nu p) = S(\nu) \left(\frac{\chi}{p}\right) p^m  \pmod{p^{2m+1}},
\]
except for very few values of $\nu$.
\begin{example} \rm
See the Ramanujan-like series \cite[eq. (1-3)]{Rama-new}. Let
\[
S(N)=\sum_{n=0}^{N-1} \frac{\left(\frac12\right)_n^5}{(1)_n^5}\left(\frac{-1}{4}\right)^n (20n^2+8n+1)
\]
If $p$ is a prime number (except for very few of them), then
\[
S(\nu p) \equiv S(\nu)  \left(\frac{1}{p}\right) p^2 \pmod{p^5},
\]
for positive integers $\nu$. Observe that for all prime $p$ we have $(1/p)=1$.

\end{example}

\begin{example} \rm
See the Ramanujan-like series \cite[eq. (4-1)]{Rama-new}. Let
\[
S(N)=\sum_{n=0}^{N-1} \frac{\left(\frac12\right)^7}{(1)_n^7}\left(\frac{1}{64}\right)^n (168n^3+76n^2+14n+1)
\]
If $p$ is a prime number (except for very few of them), then
\[
S(\nu p) \equiv S(\nu)  \left(\frac{-4}{p}\right) p^3 \pmod{p^7},
\]
for positive integers $\nu$.
\end{example}

\begin{example} \rm
See the Ramanujan-like series \cite[eq. (2-4)]{Rama-new}. Let
\[
S(N)=\sum_{n=0}^{N-1} \frac{\left(\frac12\right)_n \left(\frac12\right)_n \left(\frac13\right)_n\left(\frac23\right)_n \left(\frac16\right)_n \left(\frac56\right)_n}{(1)_n^5}\left(\frac{-1}{80^3}\right)^n (5418n^2+693n+29)
\]
If $p$ is a prime number (except for very few of them), then
\[
S(\nu p) \equiv S(\nu)  \left(\frac{5}{p}\right) p^2 \pmod{p^5},
\]
for positive integers $\nu$.
\end{example}

\section{Extended Y. Zhao's supercongruences}
By identifying numerical approximations, we conjecture that $A=r L(\chi,m+1)$, where $r$ is a rational. The $p$-adic analogue of $A$ is $A_p=r L_p(\chi,m+1)$. We have the following supercongruences:
\[
S(\nu p) \equiv \left(\frac{\chi}{p}\right) S(\nu) p^m + r z_0^{\nu} \left( \prod_{i=0}^{2m} \frac{(s_i)_{\nu}}{(1)_{\nu}} \right)  L_p(\chi,m+1)p^{2m+1} \pmod{p^{2m+2}}.
\]
which generalizes for positive integers $\nu$ the Yue Zhao's supercongruences for $\nu=1$ (author Y. Zhao at mathoverflow). To check these supercongruences use the following congruences
\[
L_p(\chi,m+1) \equiv L(\chi,2+m-p) \pmod{p}, \quad \zeta_p(m+1) \equiv \frac{{\rm bernoulli}(p-m-1)}{m+1} \pmod{p}.
\]
Observe that $L(1,m+1)=\zeta(m+1)$ and $L_p(1,m+1)=\zeta_p(m+1)$. For Bernoulli numbers associated to $\chi$ see \cite{Zagier-bernoulli}.

\begin{example} \rm
See the Ramanujan-like series \cite[eq. (1-3)]{Rama-new}. Let
\[
S(N)=\sum_{n=0}^{N-1} \frac{\left(\frac12\right)_n^5}{(1)_n^5}\left(\frac{-1}{4}\right)^n (20n^2+8n+1), \quad T(\nu)=\left(\frac{-1}{4}\right)^\nu \frac{\left(\frac12\right)_\nu^5}{(1)_\nu^5}.
\]
If $p$ is a prime number (except for very few of them), then
\[
S(\nu p) \equiv S(\nu) p^2 + 448 T(\nu) \zeta_p(3) \nu^5 p^5 \pmod{p^6},
\]
for positive integers $\nu$.
\end{example}

\begin{example} \rm
See the Ramanujan-like series \cite[eq. (4-1)]{Rama-new}. Let
\[
S(N)=\sum_{n=0}^{N-1} \frac{\left(\frac12\right)_n^7}{(1)_n^7}\left(\frac{1}{64}\right)^n (168n^3+76n^2+14n+1), \quad T(\nu)=\left( \frac{1}{64} \right)^{\nu} \frac{\left(\frac12\right)_\nu^7}{(1)_\nu^7}.
\]
If $p$ is a prime number (except for very few of them), then
\[
S(\nu p) \equiv S(\nu) \left(\frac{-4}{p}\right) p^3 + 1536 T(\nu) L_p(-4,4) \nu^7 p^7 \pmod{p^8},
\]
for positive integers $\nu$.
\end{example}

\begin{example} \rm
See the Ramanujan-like series \cite[eq. (2-4)]{Rama-new}. Let
\[
S(N)=\sum_{n=0}^{N-1} \frac{\left(\frac12\right)_n \left(\frac12\right)_n \left(\frac13\right)_n\left(\frac23\right)_n \left(\frac16\right)_n \left(\frac56\right)_n}{(1)_n^5}\left(\frac{-1}{80^3}\right)^n (5418n^2+693n+29), 
\]
and
\[
T(\nu)=\left(\frac{-1}{80^3}\right)^\nu \frac{\left(\frac12\right)_\nu \left(\frac13\right)_\nu\left(\frac23\right)_\nu \left(\frac 16\right)_\nu \left(\frac56\right)_\nu}{(1)_\nu^5}.
\]
If $p$ is a prime number (except for very few of them), then
\[
S(\nu p) \equiv S(\nu) \left(\frac{5}{p}\right)p^2 + 42000 T(\nu) L_p(5,3) \nu^5 p^5  \pmod{p^6},
\]
for positive integers $\nu$.
\end{example}

\section{An application of the extended supercongruences}  
In next examples, we will used the generalized Zudilin's supercongruences to obtain the rational parameters of the rational Ramanujan series. For that aim (except for a global rational factor) we just need taking a sufficiently large prime $p$ and $m$ values of $\nu$. In addition, we can check that there is a rational $r$ such that Zhao's supercongruences hold for that prime $p$ and those $m$ values of $\nu$. Hence $A_p=r L_p(\chi, m+1)$, and we conclude that $A=r L(\chi, m+1)$. Finally, observe that if $|z_0|>1$ then the series for $A$ is convergent.

\begin{example} \rm
We want to see that there is a series of the following form:
\[
\sum_{n=0}^{\infty}  \frac{\left( \frac12 \right)_n^5}{(1)_n^5} \frac{(-1)^n}{4^n} (a_0+a_1 n+ a_2 n^2)=t_0 \frac{\sqrt{\chi}}{\pi^2}, \quad \chi=1,
\]
where $a_0,a_1,a_2,t_0$ are positive integers. Indeed, using the Wilf-Zeilberger (WZ method) we proved that $a_0=1, a_1=8, a_2=20$. Here, from
\[
S(\nu p)-S(\nu)p^2 \equiv 0 \pmod{p^5}, \quad \nu=1,2,3, \dots,
\]
and taking $p=11$, and $\nu=1,2$, we get the linear system
\begin{align*}
103175 a_0 + 126304 a_1 + 81213 a_2& \equiv 0 \pmod{11^5}, \\
23608 a_0  + 21777 a_1 + 22319 a_2 & \equiv 0 \pmod{11^5}.
\end{align*}
Let $a_0=t$. From the above equations, we obtain
\begin{align*}
-66812987 t - 95491225 a_2 & \equiv 0 \pmod{11^4}, \\
-35044211 t - 95491225 a_1 & \equiv 0 \pmod{11^4}.
\end{align*}
Solving the equations taking into account that the inverse $\pmod{11^4}$ of $95491225$ is $12252$, we obtain
\begin{align*}
a_2=-14621 t \pmod{11^4} &= 20 t, \\ 
a_1=-14633 t \pmod{11^4} &= 8 t,
\end{align*}
Hence, the solutions are of the following form:
\[
a_0=t, \quad a_1=8t, \quad a_2=20t.  
\]
\end{example}

\begin{example} \rm
We want to know if there is a series of the following form:
\[
\sum_{n=0}^{\infty}  \frac{\left( \frac12 \right)_n  \left( \frac13 \right)_n \left( \frac23 \right)_n \left( \frac14 \right)_n \left( \frac34 \right)_n}{(1)_n^5} \frac{(-1)^n}{48^n} (a_0+a_1 n+ a_2 n^2)=t_0 \frac{\sqrt{\chi}}{\pi^2}, \quad \chi=1,
\]
and where $a_0,a_1,a_2,t_0$ are positive integers. Using the PSLQ algorithm we conjecture that $a_0=5, a_1=63, a_2=252$ and $t_0=48$. Here, from  
\[
S(\nu p)-S(\nu) p^2 \equiv 0 \pmod{p^5}, \quad \nu=1,2,3, \dots,
\]
and taking $p=13$, and $\nu=1,2$, we get the linear system
\begin{align*}
155250 a_1 + 1838 a_2 + 327490 a_0 & \equiv 0 \pmod{13^5}, \\
304350 a_1 + 329224 a_2 + 67674 a_0 & \equiv 0 \pmod{13^5}.
\end{align*}
Let $a_0=5t$. From the above equations, we obtain
\begin{align*}
26628 a_1 + 7535 t & \equiv 0 \pmod{13^4}, \\
26628 a_2 + 1579 t & \equiv 0 \pmod{13^4}.
\end{align*}
As the inverse $\pmod{13^4}$ of $26628$ is $9279$, we obtain
\begin{align*}
a_2=-28309 t \pmod{13^4} &= 252 t, \\ 
a_1=-28498 t \pmod{13^4} &= 63 t,
\end{align*}
Hence, the solutions are: $a_0=5t, \quad a_1=63t, \quad a_2=252t$.
\end{example} 

\begin{example} \rm
We want to know if there is a series of the following form:
\[
\sum_{n=0}^{\infty}  \frac{\left( \frac12 \right)_n^7}{(1)_n^7} \left(\frac{1}{64
} \right)^n (a_0+a_1 n+ a_2 n^2+a_3 n^3)=t_0 \frac{\sqrt{-\chi}}{\pi^3}, \quad \chi=-4,
\]
where $a_0,a_1,a_2,a_3,t_0$ are positive integers. Using the PSLQ algorithm, we conjecture that $a_0=1, a_1=14, a_2=76,a_3=168$ and $t_0=16$. Here, from
\[
S(\nu p) - S(\nu) \left( \frac{-4}{p} \right) p^3 \equiv 0, \pmod{p^7} \quad \nu=1,2,\dots,
\]
and taking $p=11$, and $\nu=1,2,3$, we get the equations
\[
2078533 a_1 + 9963171 a_2 + 11695266 a_3 + 16073136 a_0 \equiv 0 \pmod{11^7}, 
\]
\begin{align*}
12453192 a_1 + 988367 a_2 + 3883033 a_3 + 14086913 a_0 & \equiv 0 \pmod{11^7}, \\
17113786 a_1 + 2247378 a_2 + 4011161 a_3 + 7012796a_0 & \equiv 0 \pmod{11^7}.
\end{align*}
Let $a_0=t$. From the above equations, we obtain
\begin{align*}
7854385 a_1 + 3429250 a_2 + 19159030 t & \equiv 0 \pmod{11^4}, \\
3851936 a_1 + 8961898 a_2 + 5481146 t & \equiv 0 \pmod{11^4}.
\end{align*}
Solving the equations, we obtain
\begin{align*}
a_1 &= - 11965 t \pmod{11^4} = 14 t, \\ 
a_2 &=- 1255 t \pmod{11^4} = 76 t, \\
a_3 &=- 14473 t \pmod{11^4} = 168 t.
\end{align*}
\end{example}

\begin{example} \rm
We want to know if there is a series of the following form:
\[
\sum_{n=0}^{\infty}  \frac{\left( \frac12 \right)_n  \left( \frac13 \right)_n \left( \frac23 \right)_n \left( \frac16 \right)_n \left( \frac56 \right)_n}{(1)_n^5} \left( \frac{-1}{80^3} \right)^n (a_0+a_1 n+ a_2 n^2)=t_0 \frac{\sqrt{\chi}}{\pi^2}, \quad \chi=5,
\]
and where $a_0,a_1,a_2,t_0$ are positive integers. Using the PSLQ algorithm we conjecture that $a_0=29, a_1=693, a_2=5418$ and $t_0=128$. Here, from  
\[
S(\nu p)- S(\nu) \left( \frac{5}{p} \right) p^2 \equiv 0 \pmod{p^5}, \quad \nu=1,2,3, \dots,
\]
and taking $p=41$, $a_0=29 t$, and $\nu=1,2$, we get the linear system
\begin{align*}
76877806 a_2 + 113924268 a_1 + 43501045 t & \equiv 0 \pmod{41^5}, \\
88965067 a_2 + 84189111 a_1 + 113390736 t & \equiv 0 \pmod{41^5}.
\end{align*}
From the above equations, we obtain
\begin{align*}
38939 a_1 + 32305 t & \equiv 0 \pmod{41^3}, \\
29982 a_2 + 4321 t & \equiv 0 \pmod{41^3}.
\end{align*}
As the inverse of $38939 \pmod{41^3}$  is $55540$, we obtain
\begin{align*}
a_2=-63503 t \pmod{41^3} &= 5418 t, \\ 
a_1=-68228 t \pmod{41^3} &= 693 t,
\end{align*}
Hence, the solutions are: $a_0=29t, \quad a_1=693t, \quad a_2=5418t$.
\end{example} 

\section*{Acknowledgement}
I am very grateful to Wadim Zudilin for sharing several important ideas on the $p$ adics, and very specially for advising me to replace $x$ with $p, 2p, 3p, \dots$, and not only with $p$.

\end{document}